\documentclass[12pt]{article}
\usepackage{amssymb}
\usepackage{latexsym}
\def\part#1{\frac{\partial\phantom{q}}{\partial#1}}

\newenvironment{rmk}{\begin{trivlist}\item[]{\bf Remark:} }
{\end{trivlist}}
\newenvironment{ex}{\begin{trivlist}\item[]{\bf Example:} }
{\end{trivlist}}
\newenvironment{prf}{\begin{trivlist}\item[]{\bf Proof:} }
{\hfill $\Box$ \end{trivlist}}
\newtheorem{thm}{Theorem}
\newtheorem{prp}[thm]{Proposition}

\newcommand{\lie}[1]{\mathfrak{#1}}

\def\vol{\mathop{\rm vol}\nolimits}

\def\deg{\mathop{\rm deg}\nolimits}

\def\tr{\mathop{\rm tr}\nolimits}
\def\Ad{\mathop{\rm Ad}\nolimits}
\def\diag{\mathop{\rm diag}\nolimits}

\newcommand{\R}{\mathbf{R}}
\newcommand{\C}{\mathbf{C}}
\newcommand{\Z}{\mathbf{Z}}
\newcommand{\T}{\mathbf{T}}
\newcommand{\CP}{{\mathbf C}{\rm P}}
\begin{document}
\title{The Wess-Zumino term for a harmonic map}
\author{Nigel Hitchin\\[5pt]
\itshape  Mathematical Institute\\
\itshape 24-29 St Giles\\
\itshape Oxford OX1 3LB\\
\itshape England\\
 hitchin@maths.ox.ac.uk}
\maketitle
\begin{abstract}
\noindent We calculate the Wess-Zumino term $\Gamma(g)$ for a harmonic map $g:\Sigma\rightarrow G$ of a closed surface to a compact, simply connected, simple Lie group in terms of the energy and the holonomy of the Chern-Simons line bundle on the moduli space of flat $G$-connections. In the case $\Sigma=S^2$ we deduce that $\Gamma(g)$ is $0$ or $\pi$ and for $\Sigma=T^2$ and $G=SU(2)$ we give a formula involving hyperelliptic integrals.
\end{abstract}

\section{Introduction}
The Wess-Zumino term $\Gamma(g)$ is an invariant with values in  $\R/2\pi\Z$ associated to a smooth map $g:\Sigma\rightarrow G$ from a closed oriented  surface $\Sigma$ to a compact, simply connected, simple Lie group $G$. From the author's preferred perspective this is the holonomy around $\Sigma$ of the canonical gerbe on $G$ \cite{Hit2}, but in this paper we shall  use the more conventional definition in terms of the integral of the invariant $3$-form on $G$ over an extension of $g$ to a $3$-manifold $M$ with boundary $\Sigma$. When $G=SU(2)=S^3$ and $g$ is an embedding, $[\Gamma(g)/2\pi]$ is  the fraction of the total volume of $S^3$ enclosed by $g(\Sigma)$. Our purpose in this paper is to evaluate $\Gamma(g)$ when the map $g$ is {\it{harmonic}}. 

One motivation for this is the once-held belief \cite{Yau} that a minimal surface in $S^3$ should divide the sphere into two halves of equal volume, just like the equatorial $2$-sphere. This would be the statement that $\Gamma(g)=\pi$. The conjecture is now known to be false -- there is a counterexample in genus $11$ \cite{Kar} -- but nevertheless the questions of how one actually calculates this invariant, and what properties it has, are clearly of interest in this concrete setting. Moreover, given the knowledge we now possess on harmonic maps from the $2$-sphere or the $2$-torus, one feels that there must be an explicit formula in those cases. We shall give one.

We begin by deriving various properties of the Wess-Zumino term, without the harmonicity assumption. In particular we show that if the map $g$ takes values in a maximal torus, or  a canonically embedded symmetric space of $G$,  then $\Gamma(g)$ can only be equal to $0$ or $\pi$. 

The starting point for most work on harmonic maps of surfaces to groups is the interpretation of the equations as the vanishing of the curvature of a  family of flat $G$-connections $A({\theta})$ parametrized by the circle $e^{i\theta}$, with $A(0)$ being the trivial connection and $A(\pi)$ a connection trivialized by the map $g:\Sigma \rightarrow G$. For the case of $\Sigma=S^2$, any flat connection is trivial and we get this way a family of maps $g_{\theta}$. Here, by a simple integration, we find two results: the known integrality of $\vert \beta \vert^2 E/16\pi$ where $E$ is the energy  of the map $g$ and $\vert \beta \vert^2$ the square of the length of the highest root $\beta$ of $G$, and then the formula $\Gamma=\vert \beta \vert^2 E/16$. Thus for $S^2$, $\Gamma(g)$ is again $0$ or $\pi$. 

For a higher genus surface holonomy obstructs the trivialization of the flat connection $A({\theta})$, and we require a different approach. The most effective one for us is  classical Chern-Simons theory viewed as a topological  field theory as in \cite{Freed}. From this point of view, the one-parameter family of connections $A({\theta})$ describes a path in the moduli space ${\cal M}$ of all flat $G$-connections on the trivial bundle $\Sigma \times G$. From $\theta=0$ to $\theta=\pi$ this is a closed path $\gamma$ starting and ending at the gauge equivalence class of the trivial connection.  Our result relates the Wess-Zumino term to the holonomy of the natural connection on the Chern-Simons line bundle over ${\cal M}$. In fact we show that
$$H(\gamma) =\exp(-16i\vert \beta \vert^2E-i\Gamma).$$

The usefulness of this formula lies in the fact that we only need to know the gauge equivalence class of $A({\theta})$ to be able to calculate the holonomy. We apply this to the well-established integrable systems method of describing harmonic maps from a $2$-torus to a group $G$  and obtain for $G=SU(2)$ a formula in terms of hyperelliptic integrals. Together with the known formula for the energy, we  determine $\Gamma$. We consider some simple examples and show that in general $E$ and $\Gamma$ vary independently of each other.

The author wishes to thank D.~S.~Freed, J.~H.~Rawnsley and A.~Pazhitnov for useful discussions and EPSRC for support.

\section{The Wess-Zumino term}

Let $G$ be a compact simple, simply-connected Lie group. If $B(X,Y)$ denotes the Killing form on the Lie algebra $\lie{g}$, then the alternating form
$B(X,[Y,Z])$  
defines a bi-invariant closed $3$-form on $G$, whose de Rham class generates $H^3(G,\R)$.  A particular multiple of this gives a form $\Omega$ such that $[\Omega]$  generates $H^3(G,2\pi \Z)\cong \Z$. For $G=SU(n)$ we have
$$\Omega=\frac{1}{12\pi}\tr(g^{-1}dg)^3.$$
For a general Lie group the Lie algebra-valued $1$-form $g^{-1}dg$ is replaced  by the Maurer-Cartan form $\omega\in \Omega^1(G,\lie{g})$ and for some constant $c$  $$\Omega=cB(\omega,[\omega,\omega]).$$ 
 To determine  $c$, we note (see \cite{BS}) that the generator of $H_3(G,\Z)\cong \pi_3(G)$ can be realized by the inclusion of a distinguished subgroup $K$ isomorphic to  $SU(2)=S^3$. Its Lie algebra $\lie{k}$ is generated by the root spaces $\lie{g}^{\beta},\lie{g}^{-\beta}$ of a highest root $\beta$, and multiples of the dual element $\check \beta \in \lie{g}$. The Killing form $B$ on $\lie{g}$ restricts to a multiple $\lambda$ of the Killing form $b$ on $\lie{k}$, and so if $[\Omega]/2\pi$ generates $H^3(G,\Z)$, we must have
\begin{equation}
2\pi=\int_K \Omega =c\lambda \int_K b(\omega,[\omega,\omega])
\label{norm}
\end{equation}
For $S^3$ the Killing form is $-8g$ where $g$ is the induced metric from $\R^4$ and the $3$-form $b(\omega,[\omega,\omega])=-96 \vol_g$ 
and so
\begin{equation}
2\pi=-96 c \lambda \vol(S^3)=-192\pi^2c\lambda
\label{cl}
\end{equation}
But $\check\beta$ generates the lattice of a Cartan subalgebra of $S^3$ so $b(\check\beta,\check\beta)=-8$, and 
$$B(\check\beta,\check\beta)=\lambda b(\check\beta,\check\beta)=-8\lambda $$
and hence from (\ref{cl}) we have 
\begin{equation}
c=\frac{1}{12\pi B(\check\beta,\check\beta)}
\label{c}
\end{equation}
For $G=SU(2)$ this gives $c=-1/96\pi$.

An alternative way of writing this is to use $B$ to identify the Cartan subalgebra $\lie{h}$ with its dual, and then for each root $\alpha$,
\begin{equation}
\check \alpha=\frac{2\alpha}{B(\alpha,\alpha)}
\label{dual}
\end{equation}
The quadratic form $-B$ then defines a positive definite inner product $(\, , \,)$ on $\lie{h^*}$ and 
$$B(\check\beta,\check\beta)=-\frac{4}{\vert \beta \vert^{2}}$$
which gives
\begin{equation}
c=-\frac{1}{48\pi}\vert \beta \vert^{2} 
\label{cc}
\end{equation}
\vskip .25cm
Now let $g:\Sigma\rightarrow G$ be a smooth map of a closed oriented surface $\Sigma$ to $G$. Since $G$ is $2$-connected, if $M$ is an oriented $3$-manifold with boundary $\Sigma$, we can extend $g$ to $\tilde g: M\rightarrow G$. The Wess-Zumino term \cite{WZ} is then defined by
\begin{equation}
\Gamma(g)=\int_M \tilde g^*\Omega
\label{term}
\end{equation}
To see the dependence on the extension $\tilde g$ to $M$, suppose   $g_0$ is another such extension to $M_0$, then we form the closed $3$-manifold $N=M\cup -M_0$, and the two maps fit together to give a map $h:N\rightarrow G$. It follows that
$$\int_M \tilde g^*\Omega-\int_{M_0}  g_0^*\Omega=\int_N h^*\Omega$$
but this is of the form $2\pi n$ for an integer $n$ since the cohomology class of $\Omega$ lies in $H^3(G,2\pi \Z)$. Thus $\Gamma(g)$ is well-defined modulo $2\pi \Z$ by the map $g:\Sigma \rightarrow G$. We shall henceforth regard $\Gamma(g)$ as taking values in $\T=\R/2\pi \Z$.

\begin{ex} The simplest example is for the group $G=SU(2)=S^3$ and a surface $\Sigma$  embedded in $S^3$.  The invariant $3$-form on $S^3$ is a multiple of the volume form. The orientation on $\Sigma \subset S^3$ determines a choice of interior $M\subset S^3$ and then for the inclusion map $i:\Sigma\subset S^3$
$$\Gamma(i)=2\pi \frac{\vol(M)}{\vol(S^3)}$$
The  equatorial 2-sphere $S^2\subset S^3$ bounds a hemisphere and thus in this case $\Gamma(i)=\pi$. More generally, a  map $f:\Sigma \rightarrow S^2$ gives by composition a map $g=i\circ f$ and then 
$\Gamma(g)=(\deg f) \pi$.
\end{ex}

\begin{rmk}
In the physics literature the Wess-Zumino term is paired with the energy of the map
$$E=\frac{1}{2}\int_{\Sigma}\Vert g^{-1}dg\Vert^2$$
 in the partition function of the Wess-Zumino-Witten model
$Z=\int Dg e^{-kI(g)}$ where $I(g)$ is the complex functional 
$$I(g)=-6cE-i\Gamma$$
Note that  the energy $E$ depends on the choice of a conformal structure on $\Sigma$, but $\Gamma$ does not.
\end{rmk}
\vskip .25cm
We now prove  some basic properties of $\Gamma(g)$, concerning  maps $g$ which take values in some totally geodesic submanifolds of $G$. This is important amongst other things because a harmonic map to a totally geodesic submanifold is also harmonic as a map to the ambient space.

\begin{prp}: Let $V\subset G$ be a submanifold on which the $3$-form $\Omega$ restricts as a form to zero. Then there is a cohomology class $\gamma \in H^2(V;\T)$ such that if the image of $g:\Sigma \rightarrow G$ lies in  $V$, then $\Gamma(g)=g^*(\gamma)[\Sigma]$.
\end{prp}
\begin{prf} There is a smooth map $f:G\rightarrow G$ which retracts a tubular neighbourhood of $V$ onto $V$. Since $\Omega$ vanishes on $V$, $f^*\Omega$ vanishes in a neighbourhood of $V$. Since $f$ is the identity on $V$,  we can calculate the invariant by integrating $f^*\Omega$ over an extension of $g$. However, since $f^*\Omega$ vanishes in a neighbourhood of $V$, it defines a relative de Rham class $[\bar \Omega]\in H^3(G,V,\R)$ whose image in $H^3(G,\R)$ is $[\Omega]$. Reducing modulo $2\pi\Z$, $[\Omega]$ vanishes, so in the exact relative cohomology sequence
$$\rightarrow H^2(G,\T)\rightarrow H^2(V,\T)\stackrel \delta \rightarrow H^3(G,V,\T)\rightarrow H^3(G,\T)\rightarrow$$
we have
\begin{equation}
[\bar\Omega]=\delta(\gamma)
\label{cobound}
\end{equation}
for some $\gamma \in H^2(V,\T)$. Since $G$ is $2$-connected $H^2(G,\T)=0$ so $\gamma$ is unique.

Now extend $g:\Sigma \rightarrow V$ to a map $\tilde g$ of the $3$-manifold $M$ to $G$. In the relative exact sequence in integral homology we have an isomorphism
$$H_3(M,\Sigma;\Z)\cong H_2(\Sigma,\Z)\cong \Z$$
We obtain the Wess-Zumino term by integrating $\tilde g^*f^*\Omega$ over $M$, but this is the same as evaluating $[\bar \Omega]$ on the generator of $H_3(M,\Sigma,\Z)$ and from (\ref{cobound}), this is evaluating $g^*\gamma$ on the generator of $H_2(\Sigma,\Z)$, which is what we needed to prove.
\end{prf}

\vskip .25cm
\begin{prp} Let $T\subset G$ be a maximal torus. Let $C_{ij}$ be the Cartan matrix of $G$ and  $S_{ij}=\min (C_{ij},C_{ji})$. Then $S_{ij}$ modulo $2$ defines a skew-symmetric bilinear form on $H_1(T,\Z)$ with values in $\Z/2$, and thus an element $\gamma \in H^2(T,\Z/2)\subset H^2(T,\T)$.
If  the image of a map $g:\Sigma \rightarrow G$ lies in  $T$ then $\Gamma(g)=g^*\gamma[\Sigma]$. In particular, $\Gamma(g)=0$ or $\pi$.
\end{prp}
\begin{prf} Since $T$ is an abelian group, $[\omega,\omega]=0$ on $T$, so  in particular $B(\omega,[\omega,\omega])$ vanishes and therefore from Proposition 1 the Wess-Zumino term is given by pulling back an element $\gamma \in H^2(T,\T)$. 

 Now $H_2(T;\Z)$ is generated by the homology classes of homomorphisms from $S^1\times S^1$ to $T$, so it suffices to calculate the Wess-Zumino term  for the surface $\Sigma=S^1\times S^1$ and the map $g:S^1\times S^1\rightarrow T\subset G$ given by a homomorphism. Such a homomorphism is given by a pair of homomorphisms $h_1,h_2:S^1\rightarrow T$, and each of these is determined by a lattice point in the Cartan subalgebra $\lie{h}$ of $T$ -- the lattice $\Lambda$ generated by the simple coroots $\check\alpha_i$. 

So suppose 
$$g(e^{2\pi ix_1},e^{2\pi ix_2})=h_1h_2=\exp(x_1 \check\alpha_1)\exp (x_2 \check\alpha_2)$$
Extend $h_2:S^1\rightarrow T$ to a map $\tilde h_2$ from the disc $D^2$ into $G$ and take $M=S^1\times D^2$ with $\tilde g=h_1\tilde h_2:M\rightarrow G$. Now  $\Omega=cB(\omega,[\omega,\omega])$, and 
$$\tilde g^*\Omega=6\pi c dx_1\wedge B(\check\alpha_1,\tilde h_2^*[\omega,\omega])$$
But $d\omega+[\omega,\omega]/2=0$, so
$$\tilde g^*\Omega=-12 \pi c dx_1\wedge B(\check\alpha_1,d\tilde h_2^*\omega)=12\pi cd(dx_1\wedge B(\check\alpha_1,\tilde h_2^*\omega))$$
Using Stokes' theorem and the fact that $\tilde h_2$ restricts to the homomorphism $h_2$ on the boundary circle,
\begin{equation}
\int_{S^1\times D^2} \tilde g^*\Omega=24\pi^2 c\int_{S^1\times S^1} B(\check\alpha_1,\check\alpha_2)dx_1\wedge dx_2={2\pi}\frac{ B(\check\alpha_1,\check\alpha_2)}{B(\check\beta,\check\beta)}
\label{WZtorus}
\end{equation}
from (\ref{c}).  Now identifying $\lie{h}$ with its dual using $B$ we have 
 the Cartan matrix for simple roots $\alpha_1,\dots,\alpha_{k}$ given by  
$$C_{ij}=\frac {2(\alpha_i,\alpha_j)}{(\alpha_j,\alpha_j)}.$$
It follows from (\ref{dual}) that
$$\frac{B(\check\alpha_1,\check\alpha_2)}{B(\check\beta,\check\beta)}
=\frac{(\alpha_1,\alpha_2)}{(\alpha_1,\alpha_1)}\frac{(\beta,\beta)}{(\alpha_2,\alpha_2)}
=\frac{1}{2}C_{12}\frac{(\beta,\beta)}{(\alpha_1,\alpha_1)}$$
If all roots have the same length (as in $A_{k},D_{k},E_{k}$), this is half the Cartan matrix. More generally, the highest root $\beta$ is a long root and  $C_{12}\ne C_{21}$ iff the lengths of $\alpha_1$ and $\alpha_2$ differ, so if $\alpha_1$ is long it has the same length as $\beta$ and 
$$\frac{B(\check\alpha_1,\check\alpha_2)}{B(\check\beta,\check\beta)}=\frac{1}{2}C_{12}$$
Recalling that the off-diagonal terms of the Cartan matrix are negative, this gives in all cases
$$\frac{B(\check\alpha_i,\check\alpha_j)}{B(\check\beta,\check\beta)}=\frac{1}{2}S_{ij}=\frac{1}{2}\min (C_{ij},C_{ji})$$
and so from (\ref{WZtorus}) the Wess-Zumino term is $\pi S_{ij}$ as required.
\end{prf}
\begin{rmk} The $\Z/2$-valued bilinear form appears in the construction of a central extension of the loop group $LG$ in Chapter 4 of \cite{PS}. The central extension determines a line bundle on the loop group. In this case it is the natural  line bundle on the loop space $LM$ defined by a gerbe on $M$ \cite{Hit2}. What we are calculating here is in one language the holonomy of a flat gerbe on $T$ and in another the holonomy of a flat  line bundle on $LT$, the loops in the maximal torus.
\end{rmk}
\vskip .25cm
Another class of totally geodesic submanifolds of $G$ consists of symmetric spaces. Suppose $\sigma$ is an involution of the group $G$ and $K$ the fixed point set. A theorem of Cartan (see \cite{Hel} Theorem (8.2)) asserts that when $G$ is simply connected, $K$ is connected. Then (see \cite{CE}) the map 
$$\Phi(gK)=g\sigma(g^{-1})$$
 embeds the symmetric space $G/K$ as a totally geodesic submanifold of $G$. We shall evaluate the Wess-Zumino term for a map $g:\Sigma \rightarrow G/K\subset G$.
\begin{prp} 
 Suppose  the image of the map $g:\Sigma \rightarrow G$ lies in  the canonically embedded symmetric space $G/K$ with homology class $nx$ where $x$ is a generator of the cyclic group $H_2(G/K,\Z)$ (this is either $\Z,\Z/2$ or $0$). Then $\Gamma(g)=n\pi$ (modulo $2\pi$).  
\end{prp}
\begin{prf} The submanifold $V=G/K\subset G$ is homogeneous under the action $x\mapsto g x\sigma(g^{-1})$. Since the $3$-form $\Omega$ on $G$ is bi-invariant, we can determine its restriction to $V$ by looking at $g=e$. Here, the $\pm 1$ eigenspaces of $\sigma$ split the Lie algebra orthogonally as 
$\lie{g}=\lie{k}\oplus \lie{p}$ and
$$[\lie{k},\lie{k}]\subseteq\lie{k},\quad [\lie{k},\lie{p}]\subseteq \lie{p}, \quad [\lie{p},\lie{p}]\subseteq \lie{k}.$$
The subspace $\lie{p}$ is the tangent space of $V$ at the identity, so if $X,Y,Z \in \lie{p}$, $[Y,Z]\in [\lie{p},\lie{p}]\subseteq\lie{k}$ and by orthogonality of $\lie{k}$ and $\lie{p}$, $B(X,[Y,Z])=0$. The form $\Omega$ thus vanishes on $V$ so from  Proposition 1,  we need to find $\gamma \in H^2(V,\T)$.

Now from the fibration $G\mapsto G/K$, since $K$ is connected and $G$ simply connected, $G/K$ is simply connected and so using the Hurewicz isomorphism, $$H_2(G/K,\Z)\cong \pi_2(G/K)\cong \pi_1(K).$$
The results for the symmetric space $G/K$ can be found in \cite{BR}: $\pi_2(G/K)=\Z$ if and only if  $G/K$ is Hermitian symmetric, otherwise $\pi_2(G/K)=\Z/2$ or vanishes. 

 Since the Wess-Zumino term is in this case the evaluation of a cohomology class, it will be sufficient  to evaluate it for the map $S^2\rightarrow G$ given by a generator of $\pi_2(G/K)$. For this we again turn to \cite{BR}: take a long root $\alpha$ of the symmetric space and the copy of $SU(2)\subset G$ generated by the root spaces $\lie{g}^{\alpha},\lie{g}^{-\alpha}\subset \lie{p}$. Then $[\lie{g}^{\alpha},\lie{g}^{-\alpha}]\subset \lie{k}$ and this generates a circle subgroup $U(1)\subseteq K$. The inclusion $SU(2)\subset G$ induces a totally geodesic embedding of the $2$-sphere $SU(2)/U(1)$ into $G/K$ which generates $\pi_2(G/K)$. 

Put another way, the homotopy class is generated by the canonical embedding of the symmetric space $S^2=SU(2)/U(1)$ in $SU(2)$ followed by the inclusion homomorphism $SU(2)\subset G$. To calculate the Wess-Zumino term, we take the $3$-manifold $M$ with boundary $S^2$ to be a hemisphere in $S^3=SU(2)\subset G$. By choice, the $3$-form $\Omega$ integrates to $2\pi$ over the $SU(2)$ subgroup generated by a highest root. But the highest root is a long root thus the integral over $M$ is $\pi$.  
\end{prf}

We have seen here that in certain situations, the Wess-Zumino term takes the special values $0,\pi$ for any map, whether it is harmonic or not. As we shall see next, these values are taken also by harmonic maps of the $2$-sphere.

\section{Harmonic maps of the $2$-sphere}

Let $\Sigma$ be a closed surface with a conformal structure. Given a map $g:\Sigma \rightarrow G$ let $\alpha=g^*\omega$ be the pull-back of the Maurer-Cartan form. It is well-known (see \cite{Report}) that if we put
\begin{equation}
A({\theta})=\frac{1}{2}(1-e^{i\theta})\alpha^{1,0}+\frac{1}{2}(1-e^{-i\theta})\alpha^{0,1}
\label{family}
\end{equation}
where $\alpha=\alpha^{1,0}+\alpha^{0,1}$ and $\ast \alpha^{1,0}=i\alpha^{1,0}, \ast \alpha^{0,1}=-i\alpha^{0,1}$ then the condition for $g$ to be harmonic is that the connection $d+A(\theta)$ on the trivial bundle $\Sigma \times G$ should be flat for all $\theta$. We then have $A(0)=0$ and $A({\pi})=\alpha=g^*\omega$.

We use this approach now to calculate the Wess-Zumino term for a harmonic map $g:S^2\rightarrow G$.
\begin{thm} Let $g:S^2\rightarrow G$ be a harmonic map, with energy $E$, then
\begin{itemize}
\item
$\vert \beta \vert^2 E/8\in 2\pi \Z$
\item
$\Gamma(g)=\vert \beta \vert^2 E/16$
\end{itemize}
where $\vert \beta \vert^2$ is the length squared of the highest root $\beta$ of $G$. In particular, $\Gamma(g)=0$ or $\pi$.
\end{thm}
(The above integrality property of the energy is well-known (see \cite{Report})).
\begin{prf}
Since the $2$-sphere is simply connected, any flat bundle is trivial, so for each $\theta$ there is a map $g_{\theta}:S^2\rightarrow G$, (well-defined modulo left multiplication by a constant) such that $g^*_{\theta}\omega=A(\theta)$. Choose a smoothly varying family and consider the Wess-Zumino term $\Gamma(g_{\theta})$. By definition we  evaluate this by extending to $\tilde g_{\theta}:M\rightarrow G$ and, writing $A=\tilde g^*_{\theta}\omega$,
$$\Gamma(g_{\theta})=\int_M cB(A,[A,A])$$
Now differentiate this expression with respect to $\theta$. The flatness condition 
\begin{equation}
dA+\frac{1}{2}[A,A]=0
\label{flat}
\end{equation}
gives on differentiating 
\begin{equation}
dA'+\frac{1}{2}[A',A]+\frac{1}{2}[A,A']=0
\label{deriv}
\end{equation}
and hence
$$\frac{d}{d\theta}B(A,[A,A])=3B(A,[A,A'])=3B(A',[A,A])$$
Now from (\ref{deriv})
$$B(A,[A,A'])=-B(A,dA')$$
and from (\ref{flat}),
$$B(A,[A,A'])=B(A',[A,A])=-2B(A',dA)$$
It follows that
$$dB(A,A')=B(dA,A')-B(A,dA')=\frac{1}{2}B(A,[A,A'])$$
Hence
\begin{equation}
\frac{d}{d\theta}\Gamma(g_{\theta})=3c\int_M B(A,[A,A'])=6c\int_{S^2} B(A,A')
\label{diff}
\end{equation}
using Stokes' theorem. But from (\ref{family}),
\begin{equation}
B(A,A')=\frac{i}{2}(\cos \theta -1)B(\alpha^{1,0},\alpha^{0,1})
\label{baa}
\end{equation}
Now the energy of the map $g$ is given by
$$E=-\frac{1}{2}\int_{S^2}B(\alpha,\ast \alpha)=-i\int_{S^2}B(\alpha^{1,0},\alpha^{0,1})$$
so from (\ref{diff}) and (\ref{baa}) 
$$\frac{d}{d\theta}\Gamma(g_{\theta})=3c(1-\cos\theta)E.$$
At $\theta=0$, $g_{\theta}$  maps $S^2$ to a point, so $\Gamma(g_{0})=0$, hence integrating we obtain
\begin{equation}
\Gamma(g_{\theta})=3c(\theta-\sin\theta)E
\label{sphere}
\end{equation}
 Now  $\Gamma(g)$ is well-defined modulo $2\pi\Z$, so since $g_0=g_{2\pi}$, putting $\theta=2\pi$ in   (\ref{sphere}) we have
\begin{equation}
6c\pi E\in 2\pi \Z
\label{integer}
\end{equation}
 or equivalently from (\ref{cc}),
$$\frac{\vert \beta \vert^2 E}{16\pi}\in \Z$$
Putting $\theta=\pi$ in  (\ref{sphere}) we obtain
 $$\Gamma(g)=3c\pi E$$
 In particular, from (\ref{integer}), $\Gamma(g)=0$ or $\pi$.
\end{prf}

\section{Classical Chern-Simons theory}

If the surface $\Sigma$ has genus greater than zero, we cannot use the argument above, since the flat connection $A(\theta)$ may not be trivial -- it could have holonomy. To deal with this we have to understand the relationship between the invariant $3$-form $\Omega$ and the differential geometry of connections on surfaces and $3$-manifolds, and the appropriate setting for this is Chern-Simons theory.

Let $M$ be an oriented $3$-manifold and consider a connection on the trivial bundle $M\times G$. This is given by a Lie-algebra valued $1$-form $A$, and its curvature is $F=dA+[A,A]/2$. Define the Chern-Simons $3$-form $C(A)$ by
\begin{equation}
C(A)=6cB(A,F)-cB(A,[A,A])
\label{cs}
\end{equation}
In the case that $G=SU(n)$, we have 
$$C(A)=\frac{1}{4\pi}\tr(AdA+\frac{2}{3}A^3).$$
If $g:M\rightarrow G$ is a smooth map this describes a change of trivialization (or gauge) of the bundle and we obtain a new connection form. For $SU(n)$ it is  $\tilde A=g^{-1}Ag+g^{-1}dg$ and in general
$$\tilde A= \Ad(g)A+g^*\omega$$
where $\omega$, as before, is the Maurer-Cartan form. The Chern-Simons form is not gauge-invariant -- for $SU(n)$ we have
$$C(\tilde A)-C(A)=-\frac{1}{4\pi}\left(\frac{1}{3}\tr(g^{-1}dg)^3+d(\tr(dg g^{-1}A)\right)$$
and in general
\begin{equation}
C(\tilde A)-C(A)=-g^*\Omega-6cdB(g^*\omega,\Ad(g)A)
\label{cstilde}
\end{equation}
If $M$ is a {\it closed} $3$-manifold, it follows using  Stokes' theorem that 
$$\int_M C(\tilde A)-\int_M C(A)=-\int_M g^*\Omega \in 2\pi\Z$$
so that the integral -- the Chern-Simons invariant -- is well defined modulo $2\pi \Z$, independently of the choice of trivialization.

\begin{ex}  Let $g\rightarrow G$ be a harmonic map. We have a family of (flat) connections $A(\theta)$ on the trivial bundle, and hence we can think of $A(\theta)$ as defining a connection on the trivial bundle over the oriented closed $3$-manifold $\Sigma \times S^1$. We calculate the  Chern-Simons form using (\ref{cs}). Since $A$ has no $d\theta$ terms, $B(A,[A,A])=0$, and since $A$ is flat on $\Sigma \times \{\theta\}$, we have $F=d\theta\wedge A'$. Thus 
$$C(A)=6cB(A,A')\wedge d\theta.$$
Integrating over $\Sigma \times S^1$, using (\ref{baa}) we find
\begin{equation}
\int_{\Sigma \times S^1}C(A)=6\int_{S^1}d\theta\int_\Sigma cB(A,A')=3c\int_0^{2\pi}(1-\cos \theta)E=6\pi c E
\label{Emodz}
\end{equation}
\end{ex}
 For our application, we shall need to understand the meaning of the Chern-Simons integral for a $3$-manifold with {\it boundary}. This is described in detail in \cite{Freed} in the language of a classical topological field theory. The main features are as follows:
\begin{itemize}
\item
a connection $A$ on a principal $G$-bundle over an oriented surface $\Sigma$ defines a one-dimensional Hilbert space $Z(A)$,
\item
if $M$ is an oriented $3$-manifold with boundary $\Sigma$, the exponential of the Chern-Simons integral of an extension of the connection to $M$ defines a vector $Z(M)\in Z(A)$,
\item
if $\tilde M$ is a different $3$-manifold, then $Z(\tilde M)=e^{i\theta}Z(M)$ where $\theta$ is the Chern-Simons invariant for the closed $3$-manifold $\tilde M\cup -M$,
\item
the Hilbert spaces $Z(A(t))$ for  a connection on $\Sigma \times [0,1]$ define a line bundle $L$ over $[0,1]$ and the exponential of the Chern-Simons integral defines  the holonomy from $Z(A(0))$ to $Z(A(1))$ of a unitary connection on $L$.
\end{itemize}

\begin{rmk}
 The  Chern-Simons integral on a manifold with boundary as we have defined it is a number, but in order to  obtain this number we had to have a priori a trivialization of the bundle. The space of all connections is an infinite dimensional affine space and so we can join any connection by a  straight line from $A$ to the trivial connection. Parallel translation along this line gives a trivialization of the line bundle $L$. With respect to {\it this} trivialization, holonomy is a unit complex number and this is the number which is the exponential of the Chern-Simons integral. Without that initial choice of trivialization, the gauge transformation properties (\ref{cstilde}) basically determine the invariant interpretation  given by  the TFT viewpoint. 
\end{rmk}

Suppose the family of connections $A(t)$ on $\Sigma \times [0,1]$ are all gauge equivalent so that we have a smooth family of maps $g_t:\Sigma \rightarrow G$ with  $g_t(A(0))=A(t)$. Cap off $\Sigma \times \{1\}$ with a manifold $M_1$ with boundary $\Sigma$ to form a manifold $M=M_1\cup (\Sigma \times [0,1])$ and extend $g$ to $M$. From (\ref{cstilde}), if we integrate the Chern-Simons form over $\Sigma \times [0,1]$ we get
\begin{eqnarray}
\int_{\Sigma \times [0,1]}C(A(t))&=&-\int_Mg^*\Omega+\int_{M_{1}}g^*\Omega \nonumber\\
&+&6c\int_{\Sigma \times \{1\}}B(g^*\omega,\Ad(g)A)-6c\int_{\Sigma \times \{0\}}B(g^*\omega,\Ad(g)A)
\label{holo}
\end{eqnarray}
We saw in Section 1 that the integral of $g^*\Omega$ over a manifold with boundary the surface $\Sigma$ only depends modulo $2\pi\Z$ on $g$ restricted to the boundary. It follows from (\ref{holo}) that 
$$\exp (i\int_{\Sigma \times [0,1]}C(A(t)))$$
depends only on the connections at the end points of the interval $[0,1]$. In the TFT interpretation this means that the holonomy of the connection on the line bundle $L$ is independent of the path of gauge-equivalent connections joining them. (Again see \cite{Freed} for details). Put another way, we can say that the line bundle $L$ with connection on the space of all connections is pulled back from a bundle with connection on the moduli space ${\cal M}$ of all  connections modulo gauge equivalence. 

\section{Harmonic maps for higher genus surfaces}
Let $g:\Sigma \rightarrow G$ be a harmonic map of a surface of arbitrary genus $g$. Consider the flat connection $A(\theta)$ as in (\ref{family}) on the $3$-manifold with boundary $\Sigma \times [0,\pi]$. Here the connection $A(\theta)$ is trivial on the boundary (since $A(0)=0$ and $A(\pi)=g^{-1}dg$) but not necessarily on the interior: it may well have non-trivial holonomy there.

Choose a smooth family of maps $g_t:\Sigma \rightarrow G$ for $t\in [\pi,2\pi]$ such that  $g_\pi=g$ and $g_{2\pi}=e$, the constant map. Define a closed path $\gamma$ in the space of connections on $\Sigma \times G$ by 
$$a(\theta)=A(\theta)\quad (0\le\theta \le \pi),\quad a(\theta)=g^*_{\theta}\omega \quad (\pi \le\theta \le 2\pi).$$
The holonomy of a line bundle around a {\it closed} path is a complex number, independent of any choices. We have
\begin{thm} \label{holthm} The holonomy $H(\gamma)$ of the connection on the line bundle $L$ around the closed curve $\gamma$ is
$$H(\gamma)=\exp i(3\pi c E -\Gamma(g)).$$
\end{thm}
\begin{prf} The holonomy is the exponential of two terms: the Chern-Simons integral over $[0,\pi]$ and that over $[\pi,2\pi]$. For the first, we obtain 
  from (\ref{Emodz}) 
$$3c\int_0^{\pi}(1-\cos \theta)E=3\pi c E$$
For the second cap off $\Sigma \times \{2\pi\}$ with a manifold $M_0$ with boundary $\Sigma$ to form a manifold $M$ and extend $g_{2\pi}=e$ by the trivial map. On $M_0$ the Chern-Simons integrand is zero, so the holonomy from $L_{a(\pi)}$ to $L_{a(2\pi)}$ is represented by the Chern-Simons integral  over $M$. Now when $A=g^*\omega$, the Chern-Simons integrand (\ref{cs}) is 
$$C(A)=-cg^*B(\omega,[\omega,\omega])=-\Omega$$
so this integral is $-\Gamma(g)$ where $\Gamma(g)$ is the Wess-Zumino term of the map $g$.

The holonomy around the closed curve $\gamma$  is thus
$$H(\gamma)=\exp i (3\pi c E-\Gamma(g))$$
\end{prf}
\begin{rmk}  The holonomy we have given here is  independent of the path $g_t$, and as we argued in the last section, the identification by parallel translation along paths of gauge-equivalent connections is equivalent to producing an identification which makes the line bundle $L$ a  bundle with connection pulled back from the moduli space ${\cal M}$. Our holonomy is therefore really the holonomy of the connection on the {\it moduli space } of flat $G$-connections around the closed path defined by the gauge-equivalence class of $A(\theta)$, $0\le \theta\le \pi$. This quotient is however singular at the trivial connection, so the reader may prefer to think in terms of the holonomy around the curve in the space of all connections. The important point is that parallel translation around another loop which represents gauge-equivalent flat connections will give the same holonomy, so we can in principle replace this loop by a simpler one to find $3\pi c E-\Gamma(g)$ modulo $2\pi\Z$. We shall do this in the case of a torus in the next section.
\end{rmk}

Note that when $\Sigma=S^2$ the moduli space of flat connections is a point so $H(\gamma)$ is trivial and the theorem recovers the result of Theorem 4.

\section{Harmonic maps from a $2$-torus to $SU(2)$}

Let $\Sigma$ be a $2$-torus. First let us examine the moduli space of flat connections on $\Sigma$. The fundamental group  $\pi_1(\Sigma)=\Z^2$ which is abelian, so any flat connection has holonomy given by a pair $(h_1,h_2)$ of commuting elements in $G$, hence a pair of elements in  the maximal torus $T$. Up to gauge equivalence, a flat connection is classified by its holonomy, so the moduli space ${\cal M}$ for the torus is the quotient of $T\times T$ by the diagonal action of the Weyl group. If we represent $\Sigma=\R^2/2\pi\Z^2$ then we can write down a concrete representative for each gauge-equivalence class:
$$A=a_1dx_1+a_2dx_2$$
where $a_1,a_2\in \lie{h}$ are constants. Then the holonomy of this flat connection is $(h_1,h_2)=(\exp a_1,\exp a_2)$. 

We can calculate the holonomy of the line bundle $L$ over ${\cal M}$ by using the Chern-Simons integral for these simple representatives. Take a path 
$$A(t)=a_1(t)dx_1+a_2(t)dx_2\quad (0\le t \le 1)$$
 and consider it as a connection on $\Sigma \times [0,1]$. Then from (\ref{cs})
 $$C(A)=6cdt\wedge B(A,A')=6cdt\wedge (B(a_1,a_2')-B(a_2,a_1'))dx_1\wedge dx_2$$
 so the connection form for the line bundle $L$ is
 $$6c (B(a_1,da_2)-B(a_2,da_1))$$
 
Suppose $A(0)=0$ and $A(1)$ is gauge-equivalent to $A(0)$. This will be so if the end-point $(a_1(1),a_2(1))=(\lambda_1,\lambda_2)$ where $\lambda_1,\lambda_2$ are  lattice points in $\lie{h}$. These points  determine a homomorphism from the torus $\Sigma=\R^2/2\pi\Z^2$ to $T$ and this map $g:\Sigma \rightarrow T\subset G$ is the gauge transformation that provides the equivalence. Thus the holonomy around the closed path $\gamma$ formed by the path $A(t)$ and a path of gauge transformations to the identity is 
\begin{equation}
\exp [24\pi^2 c i \int_0^1(B(a_1,a_2')-B(a_2,a_1'))dt-i\Gamma(g)]
\label{firsthol}
\end{equation}
Now 
\begin{eqnarray}
 \int_0^1(B(a_1,a_2')+B(a_2,a_1'))dt&=& \int_0^1B(a_1,a_2)'dt\nonumber\\
&=& B(a_1(1),a_2(1))= B(\lambda_1,\lambda_2)\label{BL}
\end{eqnarray}
If $\lambda_1=\sum m_i\check\alpha_i,\lambda_2=\sum n_i\check\alpha_i$ then 
$$B(\lambda_1,\lambda_2)=  \sum_{i,j} B(\check \alpha_i,\check \alpha_j)m_in_j$$
thus from the proof of Proposition 2  
$$24 \pi^2 cB(\lambda_1,\lambda_2)=\pi\sum S_{ij}m_in_j$$
But from  Proposition 2 this is the Wess-Zumino term $\Gamma(g)$, thus using (\ref{BL}) in  (\ref{firsthol}) the holonomy can be written
\begin{equation}
  \exp (48 \pi^2 c i \int_0^1 B(a_1,a_2')dt)
  \label{finalhol}
 \end{equation} 
We can use this formula if we know explicitly the holonomy of the flat connections $A(\theta)$ given by the harmonic map. In fact, the determination of this holonomy is precisely the starting point for the treatment of harmonic maps from a torus to $SU(2)$ in \cite{Hit1}. The approach there  is to consider the holomorphic family of flat $SL(2,\C)$ connections obtained as in (\ref{family}), allowing the parameter to take values in $\C^*$ and not just the circle:
$$A({\zeta})=\frac{1}{2}(1-\zeta)\alpha^{1,0}+\frac{1}{2}(1-\zeta^{-1})\alpha^{0,1}$$
For each $\zeta$ we obtain a flat connection on the torus and its holonomy is given by commuting matrices $h_1(\zeta),h_2(\zeta)$ in $SL(2,\C)$. The eigenvalues  $(\mu_1,\mu_1^{-1})$,  $(\mu_2,\mu_2^{-1})$ of $h_1(\zeta),h_2(\zeta)$  are two-valued functions of $\zeta$ but $\mu_1,\mu_2$ become single-valued on a branched double covering of $\C^*$ which can be completed to a hyperelliptic curve $S$ double covering $\CP^1$. The eigenvalues $\mu_1,\mu_2$ are holomorphic and single valued on  $S\backslash \pi^{-1}\{0,\infty\}$, and the differentials $d\log\mu_1,d\log \mu_2$ extend to meromorphic differentials of the second kind on $S$.

If the holonomy is trivial for all $\zeta$, then the harmonic map is a  holomorphic or antiholomorphic map $f$ to an equatorial $2$-sphere. As we have seen, the Wess-Zumino term for such a map is just $(\deg f) \pi$.

In the general case, the hyperelliptic curve (the {\it spectral curve} $S$), of arithmetic genus $p$, can be described as follows. It lies in the total space of the line bundle $\pi:{\cal O}(p+1)\rightarrow \CP^1$ and is defined by an equation $\eta^2=P(\zeta)$ where $\eta$ is the tautological section of $\pi^*{\cal O}(p+1)$ over ${\cal O}(p+1)$ and $P(\zeta)$ is the pull-back of a section of ${\cal O}(2p+2)$ on $\CP^1$ (it is thus given by a polynomial of degree $2p+2$ in the affine coordinate $\zeta$). The spectral curve satisfies the following properties:
\begin{itemize}
\item
$P(\zeta)$ is a real section of ${\cal O}(2p+2)$ with respect to the real structure on $\CP^1$ induced from $\zeta \mapsto \bar \zeta^{-1}$,
\item
$P(\zeta)$ has no real zeros (i.e. no zeros on the unit circle $\zeta=\bar \zeta^{-1})$,
\item
$P(\zeta)$ has at most simple zeros at $\zeta=0,\zeta=\infty$,
\item
there exist differentials $\varphi_1,\varphi_2$ of the second kind (i.e. with zero residues) on $S$ with periods lying in $2\pi i \Z$,
\item
$\varphi_1,\varphi_2$ have double poles at $\pi^{-1}(0)$ and $\pi^{-1}(\infty)$ and satisfy $\sigma^*\varphi_i=-\varphi_i,\rho^*\varphi_i=-\bar \varphi_i$ where $\sigma$ is the involution on $S$ induced by multiplication by $-1$ in the fibres of ${\cal O}(p+1)$ and $\rho$ is the real structure on ${\cal O}(p+1)$ induced from $\zeta \mapsto \bar \zeta^{-1}$,
\item
the principal parts of $\varphi_1,\varphi_2$ are linearly independent over $\R$
\item
the holomorphic functions $\mu_1,\mu_2$ on $S\backslash \pi^{-1}\{0,\infty\}$ satisfying $\varphi_i=d\mu_i/\mu_i$ and $\mu_i\sigma^*\mu_i=1$ also satisfy the constraint $\mu_i(\xi)=1$ for $\xi \in \pi^{-1}\{1,-1\}.$
\end{itemize}
The main result  is the following:
\begin{thm} [\cite{Hit1}] Given a  line bundle of degree $(p+1)$ on a curve $S$ satisfying the conditions above, and which is quaternionic with respect to the real structure $\rho\sigma$, there is a harmonic map from a torus to $SU(2)$. Moreover every harmonic map which is not a holomorphic or antiholomorphic map to an equatorial $2$-sphere arises this way. The map is conformal if and only if $P(\zeta)$ vanishes at $\zeta=0$ and $\zeta=\infty$.
\end{thm}

We shall show how to calculate the Wess-Zumino term from the spectral curve.  Restricted to the unit circle, the functions $\mu_1,\mu_2$ define us our path in the moduli space of flat $SU(2)$ connections.
So lift the interval $\theta \in [0,\pi]\subset S^1\subset \CP^1$ to the spectral curve $S$ and then  $a_i(\theta)=\diag (\log \mu_i,-\log \mu_i)$ give us representatives for the modulus of the flat connection $A(\theta)$. We use the formula (\ref{finalhol}) for the group  $SU(2)$ and we obtain
\begin{equation}
H(\gamma)=\exp(\frac{i}{\pi} \int_0^\pi \log \mu_1  (\log \mu_2)'\, d\theta)
\label{su2hol}
\end{equation}
So from Theorem 5, we have
$$\Gamma(g)=-\frac{1}{\pi}\int_0^{\pi} \log \mu_1  (\log \mu_2)'\, d\theta-\frac{1}{32}E$$
so if we  know the energy $E$ too, we can  obtain $\Gamma(g)$. From (\ref{Emodz}), we know that
$$\frac{1}{16} E=  -\frac{1}{\pi}\int_0^{2\pi} \log \mu_1  (\log \mu_2)'\, d\theta\quad {\rm mod}\, 2\pi \Z$$ and this  will give $E/32$ modulo $\pi\Z$, but we need to know it modulo $2\pi \Z$. Fortunately there is another way of finding the energy, derived in \cite{Hit1} by using the extension of the line bundle $L$ to a holomorphic bundle over the space of $\bar \partial$ operators. We express this result in a more invariant form below,  and using $-B$ as the metric.

Recall that if $s,t$ are holomorphic sections of a line bundle $L$ on a Riemann surface, we have a well-defined {\it Wronskian} $W(s,t)=st'-ts'$ which is a holomorphic section of $L^2K$, $K$ being the canonical bundle. Take the differentials $
\varphi_1,\varphi_2$, which are meromorphic sections of $K$, and form their Wronskian $W(\varphi_1,\varphi_2)$, which is a meromorphic section of $K^3$. Each differential has a double pole at $\zeta=0$, and zero residue, so the Wronskian has a triple pole there. In a local coordinate $z$ it therefore has the form 
$$ C\frac{dz^3}{z^3}+\dots$$
The coefficient $C$, the principal part, is independent of  $z$, just like the residue of a differential. For the spectral curve, $\zeta$ is a local coordinate near $\zeta=0$ in the non-conformal case and $\eta$ is in the conformal case. The formula for the energy is then as follows:
\begin{thm} [\cite{Hit1}] Let $g$ be a  harmonic map from a $2$-torus to $SU(2)$ and $\varphi_1,\varphi_2$ the differentials of the second kind on its spectral curve. Then
\begin{itemize}
\item
if $g$ is non-conformal, the energy is given by
$$E=-16iW(\varphi_1,\varphi_2)\vert_{\zeta=0}$$
\item
if $g$ is conformal, the energy is given by
$$E=-8iW(\varphi_1,\varphi_2)\vert_{\eta=0}$$
\end{itemize}
\end{thm}

\begin{rmk} Both the formula for the energy $E$ and for the Wess-Zumino term $\Gamma$ are independent of the choice of line bundle on the spectral curve in Theorem 6. When the spectral curve is smooth, this line bundle varies on a $p$-dimensional real torus and we obtain a corresponding $p$-dimensional family of harmonic maps. The constancy of the energy  is to be expected as this is a deformation in the critical locus for the energy functional. In the conformal case, there is also a variational reason for the constancy of $\Gamma$: recall that the problem of minimizing the area of a surface {\it while fixing the enclosed volume} leads to the  equations for constant mean curvature.  Minimal surfaces, where the mean curvature vanishes, are such surfaces. 
\end{rmk}

\section{Examples}
Let us evaluate $\Gamma(g)$ for some examples, taken from \cite{Hit1}. We shall adopt the spectral curve approach, and restrict ourselves to the case $p=0$ of a rational spectral curve. These  correspond to harmonic maps of a torus to $S^3$ which are invariant under the action of a maximal torus in $SO(4)$ acting on $S^3$.

We have two formulae at our disposal, one for the energy and one for the holonomy. 
 From Theorem 7  in the nonconformal case if
$$\varphi_1=\frac{a_1}{\zeta^2}+b_1+\dots,\quad \varphi_2=\frac{a_2}{\zeta^2}+b_2+\dots$$
 then the energy is 
\begin{equation}
E=-32i(a_2b_1-a_1b_2)
\label{nonco}
\end{equation}
In the conformal case
$$\varphi_1=\frac{a_1}{\eta^2}+b_1+\dots,\quad \varphi_2=\frac{a_2}{\eta^2}+b_2+\dots$$
and 
\begin{equation}
E=-16i(a_2b_1-a_1b_2)
\label{co}
\end{equation}
From (\ref{su2hol}) and Theorem 5 we have the formula for the Wess-Zumino term:
\begin{equation}
\Gamma= -\frac{1}{\pi} \int_0^{\pi} \log \mu_1  (\log \mu_2)'\, d\theta-\frac{1}{32}E
\label{gammaformula}
\end{equation}

\subsection{The Clifford torus}
This is the simplest example of a harmonic map from a torus to $SU(2)$ -- it is an embedded minimal surface (therefore a conformal harmonic map) given by
$$g(e^{i\theta},e^{i\phi})=\frac{1}{\sqrt{2}} \pmatrix{e^{i\theta}& e^{i\phi}\cr
         -e^{-i\phi}& e^{-i\theta}}$$
         this is considered in some detail in Section 6 of \cite{Hit1}. The spectral curve  has the simple form
$$\eta^2=\zeta$$
The two functions $\mu_1,\mu_2$  have single-valued logarithms given by
$$\log \mu_1=\frac{\pi}{2}(1+i)(\eta + \frac{i}{\eta})+i\pi,\quad \log \mu_2=\frac{\pi}{2}(1-i)(\eta - \frac{i}{\eta})-i\pi
$$ and so 
$$\varphi_1=\frac{\pi}{2}(1+i)( -\frac{i}{\eta^2}+1)d\eta,\quad \varphi_2=\frac{\pi}{2}(1-i)(\frac{i}{\eta^2}+1)d\eta$$
The map is conformal so from (\ref{co})  the calculation for the energy is
$E=16\pi^2$. Now
\begin{eqnarray*}
\log\mu_1d(\log \mu_2)&=&\frac{\pi^2}{2}\frac{(\eta^2+i)^2}{\eta^3}d\eta+\frac{\pi^2}{2}(1+i)d(\eta - \frac{i}{\eta})\\
&=& \frac{\pi^2}{4}\frac{(\zeta+i)^2}{\zeta^2}d\zeta+\frac{\pi^2}{2}(1+i)d(\eta - \frac{i}{\eta})
\end{eqnarray*}
 and thus
 $$\int_0^\pi \log \mu_1  (\log \mu_2)'\, d\theta=-\frac{\pi^3}{2}-3\pi^2$$
 Thus from (\ref{gammaformula}),
 $$\Gamma=3\pi+\frac{\pi^2}{2}-\frac{\pi^2}{2}=3\pi$$
 We deduce that, modulo $2\pi$, $\Gamma=\pi$ and the Clifford torus divides the $3$-sphere into two halves of equal volume  (as can be checked by other means).

\subsection{The nonconformal rational case}
 For nonconformal harmonic maps, the spectral curve has the form
$$\eta^2=-\bar\alpha\zeta^2+(1+\alpha\bar\alpha)\zeta-\alpha$$
with $\alpha\ne 0$. This is a double covering of $\CP^1$ branched over the two points  $\alpha,\bar\alpha^{-1}$. Up to a covering (see \cite{Hit1}) a harmonic map defined  by such a curve is given by one where
$$\log \mu_1=\frac{\pi}{r}\eta\left(1-\frac{1}{\zeta}\right),\quad \log \mu_2=\frac{i \pi}{s}\eta \left(1+\frac{1}{\zeta}\right)$$
and $r=\sqrt{(1+\alpha\bar\alpha)+(\alpha+\bar \alpha)}$ and $s=\sqrt{(1+\alpha\bar\alpha)-(\alpha+\bar \alpha)}$. (One may check easily that  the constraints $\mu_i(\xi)=1$ for $\xi\in \pi^{-1}\{\pm 1\}$ are satisfied.)

Using the expansion $\eta \sim \sqrt{-\alpha}(1-(1+\alpha\bar \alpha)\zeta/2\alpha+\dots)$ near $\zeta=0$ we find the Wronskian
$$W=\frac{-i\pi^2}{rs}(1+\alpha\bar \alpha).$$
The map is non-conformal so from (\ref{nonco}) the energy calculation gives  
$$E=\frac{32\pi^2}{rs}(1+\alpha\bar \alpha).$$ In this case

$$\log\mu_1d(\log \mu_2)=\frac{i\pi^2}{rs}\eta(1-\frac{1}{\zeta})\left((1+\frac{1}{\zeta})d\eta-\eta\frac{d\zeta}{\zeta^2}\right)$$
$$=\frac{i\pi^2}{rs}\left((1-\frac{1}{\zeta^2})(-\bar\alpha \zeta+\frac{1+\alpha\bar \alpha}{2})+(\frac{1}{\zeta^2}-\frac{1}{\zeta^3})(\bar\alpha\zeta^2-(1+\alpha\bar\alpha)\zeta+\alpha)\right)d\zeta$$
using $\eta^2=-\bar\alpha\zeta^2+(1+\alpha\bar\alpha)\zeta-\alpha$ and its derivative. We now obtain the integral
$$\int_0^\pi \log \mu_1  (\log \mu_2)'\, d\theta=\frac{i\pi^2}{rs}(-i\pi(1+\alpha\bar\alpha)+2(\alpha-\bar\alpha))$$
Thus from (\ref{gammaformula}),
 $$\Gamma=-\frac{i\pi}{rs}(-i\pi(1+\alpha\bar\alpha)+2(\alpha-\bar\alpha))+\frac{\pi^2}{rs}(1+\alpha\bar \alpha)=-\frac{2\pi i}{rs}(\alpha-\bar\alpha)$$

\begin{rmk} It is interesting to calculate in this last case the WZW functional
\begin{eqnarray*}
I(g)&=&-\frac{1}{16\pi}E-i\Gamma\\
&=&-\frac{2\pi}{rs}((1+\alpha\bar \alpha)+(\alpha-\bar\alpha))\\
&=&-{2\pi}\frac{(1+\alpha\bar \alpha)+(\alpha-\bar\alpha)}{\sqrt{(1+\alpha\bar \alpha)^2-(\alpha+\bar\alpha)^2}}
\end{eqnarray*}
We see that as the branch point $\alpha$ varies in the unit disc, the real and imaginary parts of $I(g)$, which are essentially the energy and the Wess-Zumino term, vary independently. We have finally found examples where $\Gamma$ can take a continuous range of values and not just the $0,\pi$ previously encountered. 
\end{rmk}

Lastly what do we learn about the Wess-Zumino term of an {\it embedded} minimal torus? Unfortunately very little, since an old conjecture of Lawson \cite{HBL} -- that the only embedded minimal torus in $S^3$ is the Clifford torus -- remains undecided. The only obvious statement we can make is that an embedded torus cannot have $\Gamma(g)=0$ since it must enclose a positive volume. 

\end{document}